\documentstyle[12pt]{article}
\begin{document}

\title{On homeomorphisms and $C^{1}$ maps\footnote{Dedicated to the memory of my grandparents Nikolaos and Alexandra, and Konstantinos and Eleni.}}

\author{Nikolaos E. Sofronidis\footnote{$A \Sigma MA:$ 130/2543/94}}

\date{\footnotesize Department of Economics, University of Ioannina, Ioannina 45110, Greece.
(nsofron@otenet.gr, nsofron@cc.uoi.gr)}

\maketitle

\begin{abstract}
Our purpose in this article is first, following [8], to prove that
if $\alpha $, $\beta $ are any points of the open unit disc
$D(0;1)$ in the complex plane ${\bf C}$ and $r$, $s$ are any
positive real numbers such that ${\overline{D}}( \alpha ;r)
\subseteq D(0;1)$ and ${\overline{D}}( \beta ;s) \subseteq
D(0;1)$, then there exist $t \in (0,1)$ and a homeomorphism $h :
{\overline{D}}(0;1) \rightarrow {\overline{D}}(0;1)$ such that
${\overline{D}}( \alpha ;r) \subseteq D(0;t)$, ${\overline{D}}(
\beta ;s) \subseteq D(0;t)$, $h \left[ {\overline{D}}( \alpha ;r)
\right] = {\overline{D}}( \beta ;s)$ and $h = id$ on
${\overline{D}}(0;1) \setminus D(0;t)$, and second, following [9],
to prove that if $q \in {\bf N} \setminus \{ 0, 1 \} $ and ${\bf
B}({\bf 0};1)$ is the open unit ball in ${\bf R}^{q}$, while for
any $t>0$, we set $f^{(t)}( {\bf x} ) = \frac{ t {\bf x} }{ 1 +
(t-1) \Vert {\bf x} \Vert }$, whenever ${\bf x} \in {\overline{\bf
B}}({\bf 0};1)$, then $f^{(t)} \rightarrow id$ in $C^{1} \left(
{\overline{\bf B}}({\bf 0};1) , {\bf R}^{q} \right) $ as $t
\rightarrow 1^{+}$.
\end{abstract}

\section*{\footnotesize{{\bf Mathematics Subject Classification:} 03E65, 37E30, 54C35.}}

\section{Introduction}

Our purpose in this article is to prove in {\it ZF - Axiom of
Foundation + Axiom of Countable Choice} one theorem regarding
homeomorphisms of the closed unit disc ${\overline{D}}(0;1)$ in
the complex plane ${\bf C}$ and one theorem regarding maps in
$C^{1} \left( {\overline{\bf B}}({\bf 0};1) , {\bf R}^{q} \right)
$, where ${\overline{\bf B}}({\bf 0};1)$ is the the closed unit
ball in ${\bf R}^{q}$ and $q \in {\bf N} \setminus \{ 0,1 \} $.
For the definition of a homeomorphism between metric spaces see
page 39 of [1] or page 84 of [2] or page 18 of [3] or page 45 of
[4] or page 144 of [6], while the metric of $C^{1} \left(
{\overline{\bf B}}({\bf 0};1) , {\bf R}^{q} \right) $ is easily
deduced from pages 7-8 of [7] and pages 648-649 of [5], since $f =
\left( f_{1}, ..., f_{q} \right) \in C^{1} \left( {\overline{\bf
B}}({\bf 0};1) , {\bf R}^{q} \right) $, if and only if for any $i
\in \{ 1,...,q \} $, we have that $f_{i} \in C^{1} \left(
{\overline{\bf B}}({\bf 0};1) , {\bf R} \right) $.

\section{A particular homeomorphism of ${\overline{D}}(0;1)$}

{\bf 2.1. Definition.} Let $\alpha \in D(0;1)$ be arbitrary but
fixed and let $\rho > 0$ be such that ${\overline{D}}( \alpha ;
\rho ) \subseteq D(0;1)$ or (equivalently) $\vert \alpha \vert +
\rho < 1$, i.e., $0 < \rho < 1 - \vert \alpha \vert $. Finally,
let $\delta \geq 0$ be such that ${\overline{D}} \left( \alpha ;
\rho + 2 \delta \right) \subseteq D(0;1)$ or (equivalently) $\vert
\alpha \vert + \left( \rho + 2 \delta \right) < 1$, i.e., $0 \leq
\delta < \frac{ 1 - \vert \alpha \vert - \rho }{2}$. Then,
following page 5 of [8], we set
\[
{\psi }_{ \alpha ; \rho ; \delta } \left( \alpha + re^{it} \right)
= \left\{
\begin{array}{lllll}
\alpha + \frac{ \rho + \delta }{ \rho } re^{it} & \mbox{if $0 \leq
r \leq \rho $, $0 \leq t < 2 \pi $}
\\ \\
\alpha + \left( \frac{ r - \rho }{2} + \rho + \delta \right)
e^{it} & \mbox{if $\rho \leq r \leq \rho + 2 \delta $, $0 \leq t <
2 \pi $}
\\ \\
\alpha + re^{it} & \mbox{otherwise}
\end{array}
\right.
\]
whenever $\alpha + re^{it} \in {\overline{D}}(0;1)$, where $r \geq
0$ and $0 \leq t < 2 \pi $. An almost verbatim repetition of the
argument on page 5 of [8] proves that the map $$\left[ 0, \frac{ 1
- \vert \alpha \vert - \rho }{2} \right) \ni \delta \mapsto {\psi
}_{ \alpha ; \rho ; \delta } \in H \left( D(0;1) \right) $$ is
continuous, while ${\psi }_{ \alpha ; \rho ; \delta } = id$ on the
unit circle.
\\ \\
{\bf 2.2. Definition.} Let $a \in {\bf R}$ and let $0 \leq b < 1$,
while $0 \leq \epsilon < 1-b$. We set
\[
{\sigma }_{ a ; b ; \epsilon } \left( re^{i \theta } \right) =
\left\{
\begin{array}{lllll}
r \exp \left( i \left( \theta + a \right) \right) & \mbox{if $0
\leq r \leq b$, $0 \leq \theta < 2 \pi $}
\\ \\
r \exp \left( i \left( \theta - \frac{a}{\epsilon }(r-b) + a
\right) \right) & \mbox{if $b \leq r \leq b + \epsilon $, $0 \leq
\theta < 2 \pi $}
\\ \\
re^{i \theta } & \mbox{if $b + \epsilon \leq r \leq 1$, $0 \leq
\theta < 2 \pi $}
\end{array}
\right.
\]
It is not difficult to verify that ${\sigma }_{a;b; \epsilon } :
{\overline{D}}(0;1) \rightarrow {\overline{D}}(0;1)$ constitutes a
homeomorphism which is the identity on ${\overline{D}}(0;1)
\setminus D \left( 0 ; b + \epsilon \right) $ and rotation by $a$
on ${\overline{D}}( 0;b )$. Moreover, it is not difficult to
verify that if $0 < \delta < \epsilon < 1-b$, then for any $z \in
{\overline{D}}(0;1)$, we have that $\left\vert {\sigma }_{a;b;
\epsilon }(z) - {\sigma }_{a;b; \delta }(z) \right\vert < \epsilon
- \delta $.
\\ \\
{\bf 2.3. Definition.} Let $0 \leq u < 1$ be arbitrary but fixed
and let $\delta > 0$ be such that $\left[ - 2 \delta , u + 2
\delta \right] \times \left[ - 2 \delta , 2 \delta \right]
\subseteq D(0;1)$. We construct a homeomorphism $${\tau }_{u;
\delta } : {\overline{D}}(0;1) \rightarrow {\overline{D}}(0;1)$$
which is the identity on ${\overline{D}}(0,1) \setminus \left(
\left( - 2 \delta , u + 2 \delta \right) \times \left( - 2 \delta
, 2 \delta \right) \right) $ and translation by $-u$ on $\left[ u
- \delta , u + \delta \right] \times \left[ - \delta , \delta
\right] $, i.e., it translates $\left[ u - \delta , u + \delta
\right] \times \left[ - \delta , \delta \right] $ to $\left[ -
\delta , \delta \right] ^{2}$. We proceed by defining ${\tau }_{u
; \delta }$ on $\left[ - 2 \delta , u + 2 \delta \right] \times
\left[ - 2 \delta , 2 \delta \right] $. If for any $j \in \{ 1,2
\} $, we set $pr_{j} : {\bf R}^{2} \ni \left( x_{1} , x_{2}
\right) \mapsto x_{j} \in {\bf R}$, then we distinguish the
following three cases:
\begin{enumerate}
\item[(i)]
$\delta \leq y \leq 2 \delta $. If $- \sqrt{ 1 - y^{2} } \leq x
\leq \sqrt{ 1 - y^{2} }$, then we set $pr_{2} {\tau }_{u; \delta
}(x,y)=y$, while $pr_{1} {\tau }_{u; \delta }(x,y)$ is defined as
follows:
\begin{enumerate}
\item[(a)]
If $u + \delta \leq x \leq u + 2 \delta $, then
\begin{enumerate}
\item[ ]
$pr_{1} {\tau }_{u; \delta }(x,y)$
\item[ ]
$= \left( \frac{1}{ \delta } \left( 1 - \frac{ u + \delta }{
\delta } \right) ( y - \delta ) + \frac{ u + \delta }{ \delta }
\right) ( x - u - \delta ) + \frac{u}{ \delta } ( y - \delta ) +
\delta $.
\end{enumerate}
\item[(b)]
If $u - \delta \leq x \leq u + \delta $, then $pr_{1} {\tau }_{u;
\delta }(x,y) = x + \frac{u}{ \delta } ( y - \delta ) - u$.
\item[(c)]
If $- 2 \delta \leq x \leq u - \delta $, then
\begin{enumerate}
\item[ ]
$pr_{1} {\tau }_{u; \delta }(x,y)$
\item[ ]
$= \left( \frac{1}{ \delta } \left( 1 - \frac{ \delta }{ u +
\delta } \right) ( y - \delta ) + \frac{ \delta }{ u + \delta }
\right) ( x + 2 \delta ) - 2 \delta $.
\end{enumerate}
\end{enumerate}
\item[(ii)]
$- \delta \leq y \leq \delta $. If $- \sqrt{ 1 - y^{2} } \leq x
\leq \sqrt{ 1 - y^{2} }$, then we set $pr_{2} {\tau }_{u; \delta
}(x,y)=y$, while
\[
pr_{1} {\tau }_{u; \delta }(x,y) = \left\{
\begin{array}{lllll}
\frac{ u + \delta }{ \delta } ( x - u - \delta ) + \delta &
\mbox{if $u + \delta \leq x \leq u + 2 \delta $}
\\ \\
x - u & \mbox{if $u - \delta \leq x \leq u + \delta $}
\\ \\
\frac{ \delta }{ u + \delta } ( x + 2 \delta ) - 2 \delta &
\mbox{if $- 2 \delta \leq x \leq u - \delta $}
\end{array}
\right.
\]
\item[(iii)]
$- 2 \delta \leq y \leq - \delta $. If $- \sqrt{ 1 - y^{2} } \leq
x \leq \sqrt{ 1 - y^{2} }$, then we set $pr_{2} {\tau }_{u; \delta
}(x,y)=y$, while $pr_{1} {\tau }_{u; \delta }(x,y)$ is defined as
follows:
\begin{enumerate}
\item[(a)]
If $u + \delta \leq x \leq u + 2 \delta $, then
\begin{enumerate}
\item[ ]
$pr_{1} {\tau }_{u; \delta }(x,y)$
\item[ ]
$= \left( \frac{1}{ \delta } \left( 1 - \frac{ u + \delta }{
\delta } \right) ( y + 2 \delta ) + \frac{ u + \delta }{ \delta }
\right) ( x - u - \delta ) + \frac{u}{ \delta } ( y + 2 \delta ) +
\delta $.
\end{enumerate}
\item[(b)]
If $u - \delta \leq x \leq u + \delta $, then $pr_{1} {\tau }_{u;
\delta }(x,y) = x + \frac{u}{ \delta } ( y + 2 \delta ) - u$.
\item[(c)]
If $- 2 \delta \leq x \leq u - \delta $, then
\begin{enumerate}
\item[ ]
$pr_{1} {\tau }_{u; \delta }(x,y)$
\item[ ]
$= \left( \frac{1}{ \delta } \left( 1 - \frac{ \delta }{ u +
\delta } \right) ( y + 2 \delta ) + \frac{ \delta }{ u + \delta }
\right) ( x + 2 \delta ) - 2 \delta $.
\end{enumerate}
\end{enumerate}
\end{enumerate}

\noindent {\bf 2.4. Theorem.} If $\alpha $, $\beta $ are any
points of $D(0;1)$ and $r$, $s$ are any positive real numbers such
that ${\overline{D}}( \alpha ;r) \subseteq D(0;1)$ and
${\overline{D}}( \beta ;s) \subseteq D(0;1)$, then there exist $t
\in (0,1)$ and a homeomorphism $h : {\overline{D}}(0;1)
\rightarrow {\overline{D}}(0;1)$ such that ${\overline{D}}( \alpha
;r) \subseteq D(0;t)$, ${\overline{D}}( \beta ;s) \subseteq
D(0;t)$, $h \left[ {\overline{D}}( \alpha ;r) \right] =
{\overline{D}}( \beta ;s)$ and $h = id$ on ${\overline{D}}(0;1)
\setminus D(0;t)$.
\\ \\
{\bf Proof.} If $0 < \epsilon < \frac{1}{4} \min \{ \vert \alpha
\vert , 1 - \vert \alpha \vert , \vert \beta \vert , 1 - \vert
\beta \vert \} $, then it is not difficult to verify that
\[
{\overline{D}}( \alpha ; \epsilon ) = \left\{
\begin{array}{lll}
{\psi }_{ \alpha ; r ; ( \epsilon - r )/2 } \left[ {\overline{D}}(
\alpha ; r ) \right] & \mbox{if $0 < r < \epsilon $}
\\ \\
{\psi }_{ \alpha ; \epsilon ; ( r - \epsilon )/2 }^{-1} \left[
{\overline{D}}( \alpha ; r ) \right] & \mbox{if $0 < \epsilon <
r$}
\end{array}
\right.
\]
and
\[
{\overline{D}}( \beta ; \epsilon ) = \left\{
\begin{array}{lll}
{\psi }_{ \beta ; s ; ( \epsilon - s )/2 } \left[ {\overline{D}}(
\beta ; s ) \right] & \mbox{if $0 < s < \epsilon $}
\\ \\
{\psi }_{ \beta ; \epsilon ; ( s - \epsilon )/2 }^{-1} \left[
{\overline{D}}( \beta ; s ) \right] & \mbox{if $0 < \epsilon < s$}
\end{array}
\right.
\]
while if $0 \leq a < 2 \pi $ and $0 \leq b < 2 \pi $ are such that
$\alpha = \vert \alpha \vert e^{ia}$ and $\beta = \vert \beta
\vert e^{ib}$, then since ${\overline{D}}( \alpha ; \epsilon )
\subseteq {\overline{D}}( 0; \vert \alpha \vert + \epsilon )$ and
${\overline{D}}( \beta ; \epsilon ) \subseteq {\overline{D}}( 0;
\vert \beta \vert + \epsilon )$, it is not difficult to verify
that if $0 < \eta < \min \{ 1 - \vert \alpha \vert - \epsilon , 1
- \vert \beta \vert - \epsilon \} $, then ${\overline{D}}( \vert
\alpha \vert ; \epsilon ) = {\sigma }_{ - a ; \vert \alpha \vert +
\epsilon ; \eta } \left[ {\overline{D}}( \alpha ; \epsilon )
\right] $ and ${\overline{D}}( \vert \beta \vert ; \epsilon ) =
{\sigma }_{ - b ; \vert \beta \vert + \epsilon ; \eta } \left[
{\overline{D}}( \beta ; \epsilon ) \right] $. Therefore, the claim
follows from the fact that if $\epsilon $ is small enough for both
$\left[ - 2 \epsilon , \vert \alpha \vert + 2 \epsilon \right]
\times \left[ - 2 \epsilon , 2 \epsilon \right] $, $\left[ - 2
\epsilon , \vert \beta \vert + 2 \epsilon \right] \times \left[ -
2 \epsilon , 2 \epsilon \right] $ to be contained in $D(0;1)$,
then ${\tau }_{ \vert \alpha \vert ; \epsilon } \left[
{\overline{D}}( \vert \alpha \vert ; \epsilon ) \right] =
{\overline{D}}( 0 ; \epsilon ) = {\tau }_{ \vert \beta \vert ;
\epsilon } \left[ {\overline{D}}( \vert \beta \vert ; \epsilon )
\right] $. \hfill $\bigtriangleup $

\section{Particular maps in $C^{1} \left(
{\overline{\bf B}}({\bf 0};1) , {\bf R}^{q} \right) $}

{\bf 3.1. Definition.} If $t>0$, then we set
$$f^{(t)}( {\bf x} ) = \frac{ t {\bf x} }{ 1 + (t-1) \Vert {\bf x}
\Vert },$$ i.e., if $f^{(t)} = \left( f_{1}^{(t)}, ...,
f_{q}^{(t)} \right) $, then $$f^{(t)}_{i} \left( x_{1}, ..., x_{q}
\right) = \frac{tx_{i}}{ 1 + (t-1) \sqrt{ x_{1}^{2} + ... +
x_{q}^{2} } },$$ whenever $1 \leq i \leq q$ and ${\bf x} = \left(
x_{1}, ..., x_{q} \right) \in {\overline{\bf B}}({\bf 0};1)$. For
$q=2$, this function is introduced on page 1 of [9]. We remark
that $f^{(1)} = id$, while for any value of the parameter $t$, a
straightforward computation shows that $$\frac{ \partial
f_{i}^{(t)} }{ \partial x_{j} } = \frac{ t(1-t)x_{i}x_{j} }{
\sqrt{ x_{1}^{2} + ... + x_{q}^{2} } \left( 1 + (t-1) \sqrt{
x_{1}^{2} + ... + x_{q}^{2} } \right) ^{2} },$$ whenever the
indices $i$, $j$ are distinct, and $$\frac{ \partial f_{i}^{(t)}
}{ \partial x_{i} } = \frac{t}{ 1 + (t-1) \Vert {\bf x} \Vert } -
t(t-1) \frac{ x_{i}^{2} }{ \Vert {\bf x} \Vert \left( 1 + (t-1)
\Vert {\bf x} \Vert \right) ^{2} },$$ for any index $i$, where it
is not difficult to prove that the origin constitutes a removable
singularity and $\left( \frac{ \partial f_{i}^{(t)} }{
\partial x_{j} } \right) _{ {\bf x} = {\bf 0} } = 0$ and $\left(
\frac{ \partial f_{i}^{(t)} }{ \partial x_{i} } \right) _{ {\bf x}
= {\bf 0} } = t$, for any pair of distinct indices $i$, $j$.
\\ \\
{\bf 3.2. Lemma.} $\left\Vert f^{(t)} - id \right\Vert _{ \infty }
\rightarrow 0$ as $t \rightarrow 1$.
\\ \\
{\bf Proof.} Given any $t>0$ and any ${\bf x} \in {\overline{\bf
B}}({\bf 0};1)$, a straightforward computation shows that
$$f^{(t)}( {\bf x} ) - {\bf x} = \frac{ (t-1)( 1 - \Vert {\bf x}
\Vert ) {\bf x} }{ 1 + (t-1) \Vert {\bf x} \Vert }$$ and hence
$$\max\limits_{ \Vert {\bf x} \Vert \leq 1 } \left\Vert f^{(t)}(
{\bf x} ) - {\bf x} \right\Vert = \max\limits_{ \Vert {\bf x}
\Vert \leq 1 } \frac{ \vert t-1 \vert \cdot ( 1 - \Vert {\bf x}
\Vert ) \Vert {\bf x} \Vert }{ 1 + (t-1) \Vert {\bf x} \Vert } =
\vert t-1 \vert \cdot \max\limits_{0 \leq s \leq 1}
\frac{s(1-s)}{1+(t-1)s}.$$ If $t>1$, then obviously
$$\max\limits_{0 \leq s \leq 1} \frac{s(1-s)}{1+(t-1)s} \leq
\max\limits_{0 \leq s \leq 1} \left( s(1-s)\right) = \frac{1}{2}$$
and $$\left\Vert f^{(t)} - id \right\Vert _{ \infty } \leq
\frac{t-1}{2} \rightarrow 0$$ as $t \rightarrow 1^{+}$. So let us
assume that $0<t<1$. Then, a straightforward computation shows
that $$\frac{d}{ds} \left( \frac{s(1-s)}{1+(t-1)s} \right) =
\frac{ (1-t) s^{2} - 2s + 1 }{ \left( 1+(t-1)s \right) ^{2} },$$
where the roots of $(1-t) s^{2} - 2s + 1 = 0$ are $\frac{ 1 -
\sqrt{t} }{1-t}$ and $\frac{ 1 + \sqrt{t} }{1-t}$. Hence, since
$$0 < \frac{ 1 - \sqrt{t} }{1-t} < 1 < \frac{ 1 + \sqrt{t}
}{1-t},$$ it follows immediately that
\begin{enumerate}
\item[ ]
$\frac{d}{ds} \left( \frac{s(1-s)}{1+(t-1)s} \right) > 0$ on
$\left[ 0, \frac{ 1 - \sqrt{t} }{1-t} \right) $ and $\frac{d}{ds}
\left( \frac{s(1-s)}{1+(t-1)s} \right) < 0$ on $\left( \frac{ 1 -
\sqrt{t} }{1-t} , 1 \right] $,
\end{enumerate}
and consequently
$$\max\limits_{0 \leq s \leq 1}\frac{s(1-s)}{1+(t-1)s} = \frac{ \frac{ 1 -
\sqrt{t} }{1-t} \left( 1 - \frac{ 1 - \sqrt{t} }{1-t} \right) }{ 1
+ (t-1) \frac{ 1 - \sqrt{t} }{1-t} } = \frac{ \left( 1 - \sqrt{t}
\right) ^{2} }{ (1-t)^{2} },$$ which implies that $$\left\Vert
f^{(t)} - id \right\Vert _{ \infty } = (1-t) \cdot \frac{ \left( 1
- \sqrt{t} \right) ^{2} }{ (1-t)^{2} } = \frac{ 1 - 2 \sqrt{t} +
t}{1-t} \rightarrow 0$$ as $t \rightarrow 1^{-}$ and the claim
follows. \hfill $\bigtriangleup $
\\ \\
{\bf 3.3. Lemma.} For any index $i$, $\left\Vert \left(
f_{i}^{(t)} \right) _{x_{i}}' - 1 \right\Vert _{ \infty }
\rightarrow 0$ as $t \rightarrow 1^{+}$.
\\ \\
{\bf Proof.} If $t>1$, then
\begin{enumerate}
\item[ ]
$\left\Vert \left( f_{i}^{(t)} \right) _{x_{i}}' - 1 \right\Vert
_{ \infty }$
\item[ ]
$= \max\limits_{ \Vert {\bf x} \Vert \leq 1} \left\vert \frac{t}{
1 + (t-1) \Vert {\bf x} \Vert } - 1 - t(t-1) \frac{ x_{i}^{2} }{
\Vert {\bf x} \Vert \left( 1 + (t-1) \Vert {\bf x} \Vert \right)
^{2} } \right\vert $
\item[ ]
$\leq \max\limits_{ \Vert {\bf x} \Vert \leq 1} \left( \left\vert
\frac{t}{ 1 + (t-1) \Vert {\bf x} \Vert } - 1 \right\vert + t(t-1)
\frac{ x_{i}^{2} }{ \Vert {\bf x} \Vert \left( 1 + (t-1) \Vert
{\bf x} \Vert \right) ^{2} } \right) $
\item[ ]
$\leq \max\limits_{ \Vert {\bf x} \Vert \leq 1} \left\vert
\frac{t}{ 1 + (t-1) \Vert {\bf x} \Vert } - 1 \right\vert + t(t-1)
\max\limits_{ \Vert {\bf x} \Vert \leq 1 } \frac{ x_{i}^{2} }{
\Vert {\bf x} \Vert \left( 1 + (t-1) \Vert {\bf x} \Vert \right)
^{2} }$
\item[ ]
$\leq \max\limits_{ \Vert {\bf x} \Vert \leq 1} \left\vert
\frac{t}{ 1 + (t-1) \Vert {\bf x} \Vert } - 1 \right\vert + t(t-1)
\max\limits_{ \Vert {\bf x} \Vert \leq 1 } \frac{ \Vert {\bf x}
\Vert ^{2} }{ \Vert {\bf x} \Vert \left( 1 + (t-1) \Vert {\bf x}
\Vert \right) ^{2} }$
\item[ ]
$\leq \max\limits_{ \Vert {\bf x} \Vert \leq 1} \left\vert
\frac{t}{ 1 + (t-1) \Vert {\bf x} \Vert } - 1 \right\vert + t(t-1)
\max\limits_{ \Vert {\bf x} \Vert \leq 1 } \frac{ \Vert {\bf x}
\Vert }{ \left( 1 + (t-1) \Vert {\bf x} \Vert \right) ^{2} }$
\item[ ]
$\leq (t-1) + t(t-1) \max\limits_{ \Vert {\bf x} \Vert \leq 1 }
\Vert {\bf x} \Vert $
\item[ ]
$\leq (t-1) + t(t-1)$
\item[ ]
$= t^{2} - 1$
\end{enumerate}
and the claim follows. \hfill $\bigtriangleup $
\\ \\
{\bf 3.4. Lemma.} For any pair of distinct indices $i$, $j$, we
have that
\begin{enumerate}
\item[ ]
$\left\Vert \left( f_{i}^{(t)} \right) _{x_{j}}' \right\Vert _{
\infty } \rightarrow 0$ as $t \rightarrow 1^{+}$.
\end{enumerate}

\noindent {\bf Proof.} If $t>1$, then
\begin{enumerate}
\item[ ]
$\left\Vert \left( f_{i}^{(t)} \right) _{x_{j}}' \right\Vert _{
\infty }$
\item[ ]
$= \max\limits_{ \Vert {\bf x} \Vert \leq 1} \left\vert \frac{
t(1-t)x_{i}x_{j} }{ \sqrt{ x_{1}^{2} + ... + x_{q}^{2} } \left( 1
+ (t-1) \sqrt{ x_{1}^{2} + ... + x_{q}^{2} } \right) ^{2} }
\right\vert $
\item[ ]
$= \max\limits_{ \Vert {\bf x} \Vert \leq 1} \left\vert \frac{
t(1-t)x_{i}x_{j} }{ \Vert {\bf x} \Vert \left( 1 + (t-1) \Vert
{\bf x} \Vert \right) ^{2} } \right\vert $
\item[ ]
$= t(t-1) \max\limits_{ \Vert {\bf x} \Vert \leq 1} \frac{
\left\vert x_{i} \right\vert \cdot \left\vert x_{j} \right\vert }{
\Vert {\bf x} \Vert \left( 1 + (t-1) \Vert {\bf x} \Vert \right)
^{2} }$
\item[ ]
$= t(t-1) \max\limits_{ \Vert {\bf x} \Vert \leq 1} \frac{ \Vert
{\bf x} \Vert \cdot \Vert {\bf x} \Vert }{ \Vert {\bf x} \Vert
\left( 1 + (t-1) \Vert {\bf x} \Vert \right) ^{2} }$
\item[ ]
$= t(t-1) \max\limits_{ \Vert {\bf x} \Vert \leq 1} \frac{ \Vert
{\bf x} \Vert }{ \left( 1 + (t-1) \Vert {\bf x} \Vert \right) ^{2}
}$
\item[ ]
$= t(t-1) \max\limits_{ \Vert {\bf x} \Vert \leq 1} \Vert {\bf x}
\Vert $
\item[ ]
$\leq t(t-1)$
\end{enumerate}
and the claim follows. \hfill $\bigtriangleup $
\\ \\
{\bf 3.5. Theorem.} $f^{(t)} \rightarrow id$ in $C^{1} \left(
{\overline{\bf B}}({\bf 0};1) , {\bf R}^{q} \right) $ as $t
\rightarrow 1^{+}$.
\\ \\
{\bf Proof.} It is an immediate consequence of the previous three
lemmas. \hfill $\bigtriangleup $

\end{document}